\documentclass[12pt]{article}
\usepackage{fullpage}
\usepackage{amsmath}
\usepackage{amssymb}
\usepackage{cite}
\usepackage{graphicx}
\usepackage{appendix}
\def\half{{\textstyle {1\over 2}}}

\def\ds{\displaystyle}

\def\R{\mathbb{R}}
\def\barr{\begin{array}}
\def\earr{\end{array}}

 \def\beq{\begin{equation}}
 \def\eeq{\end{equation}}

\newtheorem{thm}{Theorem}[section]
\newtheorem{Lemma}[thm]{Lemma}
 \newcommand{\proofend}{\vrule height 6pt width 6pt depth 1pt}
  \newcommand{\Label}{\label}
 
 \numberwithin{equation}{section}
\begin{document}
\title{Traveling Waves for  Conservation Laws with Cubic Nonlinearity and BBM Type Dispersion }
 \author{Michael Shearer\footnote{Corresponding author. Dept. of
    Mathematics,
    N.C. State University, Raleigh, NC 27695. Tel. 919-515-3298. email: shearer@ncsu.edu}, \ Kimberly R. Spayd\footnote{Department of Mathematics, Gettysburg College, Gettysburg, PA 17325} \ and \ Ellen R. Swanson\footnote{Department of Mathematics, Centre College, Danville, KY 40422} }
\maketitle
\begin{abstract}
  Scalar conservation laws with non-convex fluxes have shock wave solutions that violate the Lax entropy condition. In this paper, such solutions are selected by showing that some of them have corresponding traveling waves for the equation supplemented with dissipative and dispersive higher-order terms. For a cubic flux, traveling waves can be calculated explicitly for linear dissipative and dispersive terms.  Information about their existence can be used to solve the Riemann problem, in which we find solutions for some data that are different from the classical Lax-Oleinik construction. We consider dispersive terms of a BBM type and show that the calculation of  traveling waves is somewhat more intricate than for a KdV-type dispersion. The explicit calculation is based upon  the calculation of parabolic invariant manifolds for the associated ODE describing traveling waves. The results extend to the p-system of one-dimensional elasticity with a cubic stress-strain law.

  \end{abstract}
\section{Introduction}\label{intro}
There has been   interest recently \cite{fan, gray1, spayd} in models related to the Buckley-Leverett equation \cite{buckley} of two-phase porous media flow, in which a rate-dependent dispersive term is included in the capillary pressure.   The Buckley-Leverett flux is a non-convex fractional flow rate, so that with both dissipation and dispersion, there are likely to be undercompressive shocks \cite{hayes, lefloch}. 

In this paper, we consider a simpler equation that has a cubic flux, and both dissipation and rate dependent dispersion, similar to the BBM (Benjamin-Bona-Mahoney) equation \cite{benjamin1}. Such equations have been termed pseudo-parabolic \cite{cuesta}, and the specific equation considered we refer to as the modified BBM-Burgers equation, meaning that the BBM equation (which has a quadratic flux)  is modified here to have a  cubic flux function, and the dissipation is of the simple form seen in Burgers' equation. 

In previous work on undercompressive shocks \cite{shearer1995}, we characterized traveling wave solutions for the modified KdV-Burgers equation with explicit formulas. Here, we use  a similar analysis to identify invariant parabolic curves through equilibria for the vector field whose heteroclinic orbits represent traveling waves. For the BBM-type dispersion however, the analysis leads to an implicit parameterization, but again with explicit formulas. Consequently, the proof of uniqueness, while relying on a general result for the type of cubic vector field that arises here, is slightly less direct. 

Interestingly, there is a trade-off between the existence of undercompressive traveling waves and the stability or instability  of constant solutions that does not arise in the case of KdV-type dispersion. Consequently the cubic nonlinearity has to be chosen carefully in order that constant solutions are stable for the PDE, while preserving the existence of traveling waves approximating undercompressive shocks.

 We solve the Riemann problem explicitly for the underlying conservation law, using traveling wave information from the full equations to determine admissible shocks. 
Numerical simulations verify that these solutions give the structure of smooth solutions when dissipation and dispersion are included.
 
 In Section~\ref{sec_bbm} we describe basic properties of the PDE. In Section~\ref{sec_tw} we calculate invariant manifolds on which the  traveling  waves exist.  This information is used to solve Riemann problems in Section~\ref{sec_rp}.
In Section \ref{sec_system}  we show how a similar analysis applies to the quasilinear wave equation of one-dimensional elasticity with cubic stress-strain law including rate dependence (viscosity) and capillarity.  In contrast to the situation for non-monotonic stress-strain laws, Riemann problem solutions cannot include more than one undercompressive shock.

 \section{The Modified  BBM-Burgers Equation}\label{sec_bbm}
The BBM equation \cite{benjamin1}
\beq\label{bbm1}
u_t+uu_x+\mu u_{xxt}=0,
\eeq
is a variation on the KdV equation for water waves in a long -wave approximation. 
A similar dispersive regularization of the Buckley-Leverett equation of two-phase flow in porous media was introduced by Gray and Hassanizadeh \cite{gray1,spayd}:
\beq\label{bl1}
u_t+f(u)_x=\beta (h(u)u_x)_x +\mu (h(u)u_{xt})_x.
\eeq
In this equation, $u(x,t)$ represents a   fraction of the local pore volume occupied by of one of the two fluid phases and the flux function $f(u)$  is the fractional flow rate of that fluid, derived from Darcy's law. 
Significantly, $f(u)$ is non-convex.
The two terms on the right hand side of   equation \eqref{bl1}    are dissipative and dispersive, respectively, with positive parameters $\beta, \mu,$ and a positive nonlinear  function $h(u).$ They represent the equilibrium and rate-dependent contributions of interfacial energy to capillary pressure.

Burgers' equation 
\beq\label{burgers1}
u_t+uu_x=\beta u_{xx}
\eeq
is dissipative, and is a prototype of conservation laws regularized by viscosity $\beta>0.$ 
In our analysis of undercompressive shocks for scalar equations, we consider a greatly simplified version of equation \eqref{bl1} in which   the dissipative and dispersive terms are linear:
\beq
\label{pde2}
u_t+f(u)_x=\beta u_{xx}+\mu u_{xxt},
\eeq
and the flux function $f(u)$ is cubic. 

This equation has two crucial features when $\beta>0$ and $\mu>0.$ First, we observe that 
 \eqref{pde2} is linearly stable at every constant $u=\overline{u}.$ Specifically, the dispersion relation is
$$
\lambda +if'(\overline{u}) \xi=-\beta \xi^2-\mu \xi^2\lambda.
$$
Thus,
$$
\lambda(\xi)=\frac{-i\xi f'(\overline{u}) -\beta \xi^2}{1+\mu \xi^2}.
$$
Specifically, for $\mu>0,$ we have $Re \, \lambda<0$ for all $\xi>0.$  
However, if $\mu<0,$ then the equation is unstable, having $Re \, \lambda>0$ for all $\xi> 1/\sqrt{|\mu|}.$ 

The second property we require is the presence of  traveling waves with positive speed corresponding  to heteroclinic orbits between saddle-point equilibria of an ODE system. With $\mu>0,$ this property leads to the choice of cubic flux function
\beq\label{flux1}
f(u)=u-u^3.
\eeq
This flux function resembles the Buckley-Leverett flux in the interval $-1/\sqrt{3}<u< 1/\sqrt{3}$ in that $u-u^3$ is monotonically increasing, concave for $u<0$ and convex for $u>0.$ 
Equations such as \eqref{bl1},\eqref{pde2} are sometimes referred to as  {\em pseudoparabolic} \cite{fan,cuesta}, based on properties of the dispersion relation of the equation linearized about a constant. 

The parameter $\beta$ is non-negative in both \eqref{bl1} and \eqref{pde2} in order that the dissipative term is not destabilizing. 
The transformation $x\to -x$ changes only the sign of the flux, thereby switching both the convexity and the monotonicity, but leaving the sign of $\beta $ and $\mu$ unchanged. With this transformation, waves propagate to the left (with negative speed)  rather than to the right.  
The behavior is different for KdV-type equations, where the transformation $x\to -x$ changes the sign of both the flux and the dispersive term, for example in the modified KdV-Burgers equation \cite{shearer1995},
$$
u_t+(u^3)_x=\beta u_{xx}-\mu u_{xxx}.
$$
In summary, the equation we consider here, 
\beq
\label{pde3}
u_t+(u-u^3)_x=\beta u_{xx}+\mu u_{xxt},
\eeq
is referred to in this paper as the modified BBM-Burgers equation, the term {\em modified} being used because the quadratic flux of the BBM-Burgers equation has been replace by a cubic.
 \section{Traveling Waves}\label{sec_tw}
 The scalar conservation law
\beq
\Label{cl}
u_t+(u-u^3)_x=0
\eeq
has concave-convex flux $f(u)=u-u^3.$ The characteristic speed is $f'(u)=1-3u^2,$ so that rarefaction waves $u(x,t)=\overline{u}(x/t)$ centered at $x=t=0$   are  given by 
$\overline{u}(r)=\pm \sqrt{(1-r)/3}, \ r=x/t<1.$   

A {\em shock wave} from $u_-$ to $u_+$ with speed $s$ is a discontinuous weak solution of the scalar conservation law 
 which has the form
\beq
\Label{shock1}
u(x,t)=\begin{cases}
u_-&\text{if \ $x<st$}\\[8pt]
u_+&\text{if \ $x>st$}.
\end{cases}
\eeq
The shock speed  $s$ is defined by the Rankine-Hugoniot condition 
\beq
\Label{rh}
-s(u_+ - u_-)+u_+-u_+^3-(u_--u_-^3)=0, \quad \mbox{so that} \ \  s=1-(u_+^2+u_+u_-+u_-^2)
\eeq
is the slope of the chord connecting $(u_-,u_-- u_-^3)$ and $(u_+,u_+- u_+^3)$.  (In this paper we consider only constant $u_\pm$ and $s$;  more generally, $u_\pm$ would be one-sided limits at a discontinuity $x=\tilde x(t)$ with speed $s(t)=\tilde x^\prime(t).$)

A shock wave satisfies the Lax entropy condition if characteristics approach the shock from both sides:
$$
1-3u_+^2<s<1-3u_-^2.
$$ 
In this case, the shock is referred to as {\em compressive,} or as a {\em Lax shock.}

We shall say a shock wave is {\em TW-admissible} if there is a traveling wave solution 
\beq\label{TW1}
u(x,t)=\tilde u(\eta), \quad  \eta=x-st
\eeq
 of \eqref{pde3} that satisfies far-field conditions
\beq\label{bcs1}
\tilde{u}(\pm\infty)= u_\pm, \quad \tilde{u}^\prime(\pm\infty)=0, \quad 
\tilde{u}^{\prime\prime}(\pm\infty)=0.
\eeq
Lax shocks with $|u_--u_+|$ small (weak shocks) are TW-admissible, but stronger shocks need not be. A shock for which characteristics pass through the shock necessarily fails the Lax entropy condition. Such shocks that are TW-admissible are called {\em undercompressive.} Because the flux function $f(u)=u-u^3$ has a single inflection point, undercompressive shocks correspond to a chord cutting the graph of $f,$ as $u$ goes from $u_-$ to $u_+.$ Then the shock speed is necessarily greater than the characteristic speed on each side, so that characteristics pass through from ahead of the shock to behind. In the language of gas dynamics, the shock is supersonic with respect to the sound speed both ahead of and behind the shock.

 Substituting (\ref{TW1}) into \eqref{pde3} gives the third order ODE (omitting tildes)
\beq
\Label{tw1}
-s u^\prime+(u-u^3)^\prime=\beta u^{\prime\prime}-\mu su^{\prime\prime\prime}
\eeq
where $ ^\prime=d/d\eta.$ 
Integrating \eqref{tw1} and implementing the boundary condition at $\eta=-\infty$ yields the second order ODE 
 \beq
\Label{tw2}
-s(u-u_-)+u-u^3-(u_--u_-^3)=\beta u^\prime-\mu su^{\prime\prime}.
\eeq
It is convenient to rescale $\eta$ to eliminate the parameter $\mu s$ in the final term. For $s>0,$ we set $\xi=\eta/\sqrt{\mu s}$. Then the ODE \eqref{tw2} becomes an equation for $u=u(\xi):$
 \beq
\Label{tw3}
u^{\prime\prime}=\frac{\beta}{\sqrt{\mu s}}u^\prime+u^3-u-(u_-^3-u_-)+s(u-u_-),
\eeq
where $^\prime=d/d\xi$.
We analyze  \eqref{tw3} as a first order autonomous system with parameters $s, u_-,$ together with the combined parameter $\gamma=\frac{\beta}{\sqrt{\mu}}:$
\begin{subequations}\Label{ode1}
\begin{eqnarray}
\Label{ode1a}
u^\prime&=&v\\
\Label{ode1b}
v^\prime&=&\frac{\gamma}{\sqrt{s}} v+u^3-u-(u_-^3-u_-)+s(u-u_-).
\end{eqnarray}
\end{subequations}

Equilibria for this system \eqref{ode1} are points $(u,v)=(u,0),$ where  $u^3-u-(u_-^3-u_-)+s(u-u_-)=0$;  these correspond to points of intersection between the graph   of $f(u)=u-u^3$ and the line with slope $s$ through  $(u_-,u_--u_-^3).$   Thus, either $u=u_-$ with $s$ arbitrary, corresponding to a constant solution, or $u=u_+$ satisfies the Rankine-Hugoniot condition (\ref{rh}). Solving for $u_+,$ we find 
\beq\label{equil2}
u_+=\half\left\{-u_-\pm\sqrt{4(1-s)-3u_-^2}\right\}.
\eeq
Thus, there are exactly three equilibria when 
\beq\label{three}
1-s>3u_-^2/4.
\eeq  
{\bf Remark:} \ 
  Returning to equation \eqref{tw1}, we note that when there are three equilibria, the outside equilibria are saddle points if and only if $s> 0.$ To see this, we calculate the eigenvalues at an equilibrium. Let $c(u)=u^3-u-(u_-^3-u_-)+s(u-u_-).$ Then the eigenvalues at an equilibrium $u$ are 
  $$
  \lambda_\pm=\frac 12\left\{\frac{\beta}{\mu s} \pm\sqrt{\left(\frac{\beta}{\mu s} \right)^2+4\frac{c'(u)}{\mu s}}\right\}.
  $$
  At the outside equilibria, $c'(u)>0,$ i.e., $s>1-3u^2$.  Consequently, the eigenvalues are of opposite sign only if $s>0.$ For $s=0, $ the equation degenerates to a first order ODE. \\

Let $u_->0.$ 
In terms of $u_+,$ there are three equilibria $u_->u_0>u_+$  when 
 \beq\label{three2}
-2u_-<u_+<-u_-/2, \quad \mbox{if} \ \ 0<u_-<1/\sqrt{3},
\eeq
  and 
 \beq\label{three3}
-\half\left(u_- +\sqrt{1-3u_-^2/4}\right)<u_+<-u_-/2 \quad \mbox{if} \ \ 1/\sqrt{3}<u_-<2/\sqrt{3}. \eeq
 Moreover, we then have $1-3u_\pm^2<s<1-3u_0^2.$ Writing 
\beq\label{eq3a}
u^3-u-(u_-^3-u_-)+s(u-u_-)=(u-u_-)(u-u_0)(u-u_+),
\eeq
 we observe that there is no quadratic term on the left side, so that
\beq\label{eq4a}
u_-+u_0+u_+=0.
\eeq
 The Jacobian matrix of the vector field on the right hand side of  \eqref{ode1} is 
\beq
\Label{jac}
J=
\begin{pmatrix}
0&1\\[10pt]
s-1+3u^2 & \displaystyle \frac{\gamma}{\sqrt{s}}
\end{pmatrix}
\eeq
with eigenvalues 
\beq
\Label{eig}
\lambda_{\pm}=\frac{1}{2}\left[\frac{\gamma}{\sqrt{s}}\pm\sqrt{\frac{\gamma^2}{s}+4(s-1+3u^2)}\right].
\eeq
As observed in the above remark, the outside equilibria $u=u_\pm$ are saddles since $\lambda_\pm$ are real and of opposite sign ($3u_\pm^2>1-s$ when $u_+<u_0<u_-$ ).

\subsection{Saddle-Saddle Connections}

By a  {\em saddle-saddle connection} from $u_-$ to $u_+$ we mean a heteroclinic orbit from $(u_-,0)$ to $(u_+,0)$ when $(u_\pm,0)$ are saddle point equilibria.  
We seek  saddle-saddle connections between equilibria $u=u_\pm, v=0$ that lie on an invariant parabola
\beq
\Label{eq2}
v=k(u-u_-)(u-u_+).
\eeq
We shall find an equation relating $u_+, u_-$ in order for (\ref{eq2}) to be invariant, and $k$ will be  determined. However, 
since $v=u',$ we must have $k>0$ if $u_->0>u_+,$ guaranteeing that $u(\xi)$ decreases from $u_-$ to $u_+.$ Similarly,  $k<0$ if $u_+>0>u_-.$  

Since a saddle-saddle trajectory is necessarily a graph $v=v(u),$ we can rewrite system
(\ref{ode1}) as a single equation
\beq\label{vueqn}
v\frac{dv}{du}=\frac{\gamma}{\sqrt{s}}v+u^3-u-(u_-^3-u_-)+s(u-u_-).\\
\eeq
Now substitute (\ref{eq3a}), (\ref{eq2}) into this equation to obtain (after canceling factors $(u-u_-)(u-u_+)$):
\beq
k^2(2u-u_--u_+)=\frac{\gamma k}{\sqrt{s}}+(u-u_0).
\eeq
Consequently,
$$
2k^2=1 \quad \mbox{and} \quad  -k^2(u_-+u_+)=\frac{\gamma k}{\sqrt{s}}-u_0.
$$
Thus, $k=1/\sqrt{2},$ since $u_->0,$ so that $v(u)<0$ between $u_+$ and $u_-.$  Using (\ref{eq4a}) we find     
\beq
\Label{uzero}
u_0=\frac{\sqrt{2}}{3\sqrt{s}}\gamma>0.
\eeq
 Using (\ref{eq4a}) again and  $s=1-(u_+^2+u_-u_++u_-^2),$ we have the equation relating $u_\pm:$ 
\beq
\Label{eq6}
\sqrt{1-(u_+^2+u_-u_++u_-^2)}\,\,(u_++u_-)=-\frac{\sqrt{2}}{3}\gamma
\eeq
 We solve this equation parametrically.
   Let $u_-=-au_+,$ with $a>0$ a parameter. Then the restriction (\ref{three2}) corresponds to $\half<a<2.$ However, $a$ is further restricted by the following calculation. Suppose there is a trajectory $v=v(u)\leq 0, \ u_+\leq u\leq u_-,$ with $v(u_\pm)=0.$ 
 Integrating (\ref{vueqn}) from $u_+$ to $u_-$ yields,
 $$
 0=\frac{\gamma}{\sqrt{s}}\int_{u_+}^{u_-}
 v(u) \, du +\int_{u_+}^{u_-} (u^3-u-(u_-^3-u_-)+s(u-u_-))\, du.
 $$
 In this equation,   $v(u)=u'(\xi)<0,$ so that the first term on the right is negative. Consequently,
 \beq\label{integral1}
 \int_{u_+}^{u_-} (u^3-u-(u_-^3-u_-)+s(u-u_-))\, du>0.
 \eeq
 This simply means the signed area between the curve $y=u-u^3$ and the chord 
 $y=s(u-u_-)+u_- - u_-^3$ is negative. This area is zero precisely when $u_+=-u_-, $
 and remains negative in the interval $-2u_-<u_+<-u_-,$ or $-u_+/2<u_-<u_+.$  Thus, we must have 
 \beq\label{arest}
 \displaystyle \frac{1}{2}<a<1.
 \eeq
  Substituting $u_-=-au_+$ into equation (\ref{eq6}), we obtain 
  \beq\label{eq7}
 u_+(1-a) \sqrt{1-(a^2-a+1)u_+^2}=-\frac{\sqrt{2}}{3}\gamma
  \eeq
Let $X=u_+^2.$ Then $u_+=-\sqrt{X},$ so that
  \beq\label{eq8}
 \sqrt{X}(1-a) \sqrt{1-(a^2-a+1)X}=\frac{\sqrt{2}}{3}\gamma
\eeq
 Squaring both sides, we obtain a quadratic equation for $X,$ with solutions
 $$
 X=\ds\frac{1\pm\sqrt{D(a,\gamma)}}{2(1-a+a^2)}, \quad D(a,\gamma)=1-\frac89\gamma^2\left(1+\frac{a}{(a-1)^2}\right).
 $$
 In this expression, we require $D(a,\gamma)\geq 0,$ and we obtain   
  a parametric representation of the curve for   $\gamma>0,$
   \beq\label{parametric1}
 u_+=u_+(a)=-\left(\ds\frac{1\pm\sqrt{D(a,\gamma)}}{2(1-a+a^2)}\right)^\half,   \ u_-=-au_+(a). 
   \eeq
$D(a,\gamma)$ is positive for at most part of the range $\half<a<1.$ The following lemma is proved directly by investigating the zeroes of $D(a,\gamma).$  
 
 \begin{Lemma} Let $k=\ds\frac{8\gamma^2}{9}$ and $\tilde{a}= \ds\frac{1}{2(k-1)}(k-2+\sqrt{k(4-3k)}).$ \\[6pt]
 (i) \ If $0<\gamma<\sqrt{3/8},$ then $1/2<\tilde{a}<1, \ 0<D(a,\gamma)<1-\frac 83 \gamma^2$ for $\half<a<\tilde{a},$    and $D(a,\gamma)<0$ for $\tilde{a}<a<1.$\\[6pt]
 (ii) \ If $\gamma>\sqrt{3/8},$ then $D(a,\gamma)<0$ for all $a\in [\half,1].$ 
  \end{Lemma}

\begin{figure}[ht]
\centerline{\includegraphics[height=3.7in]{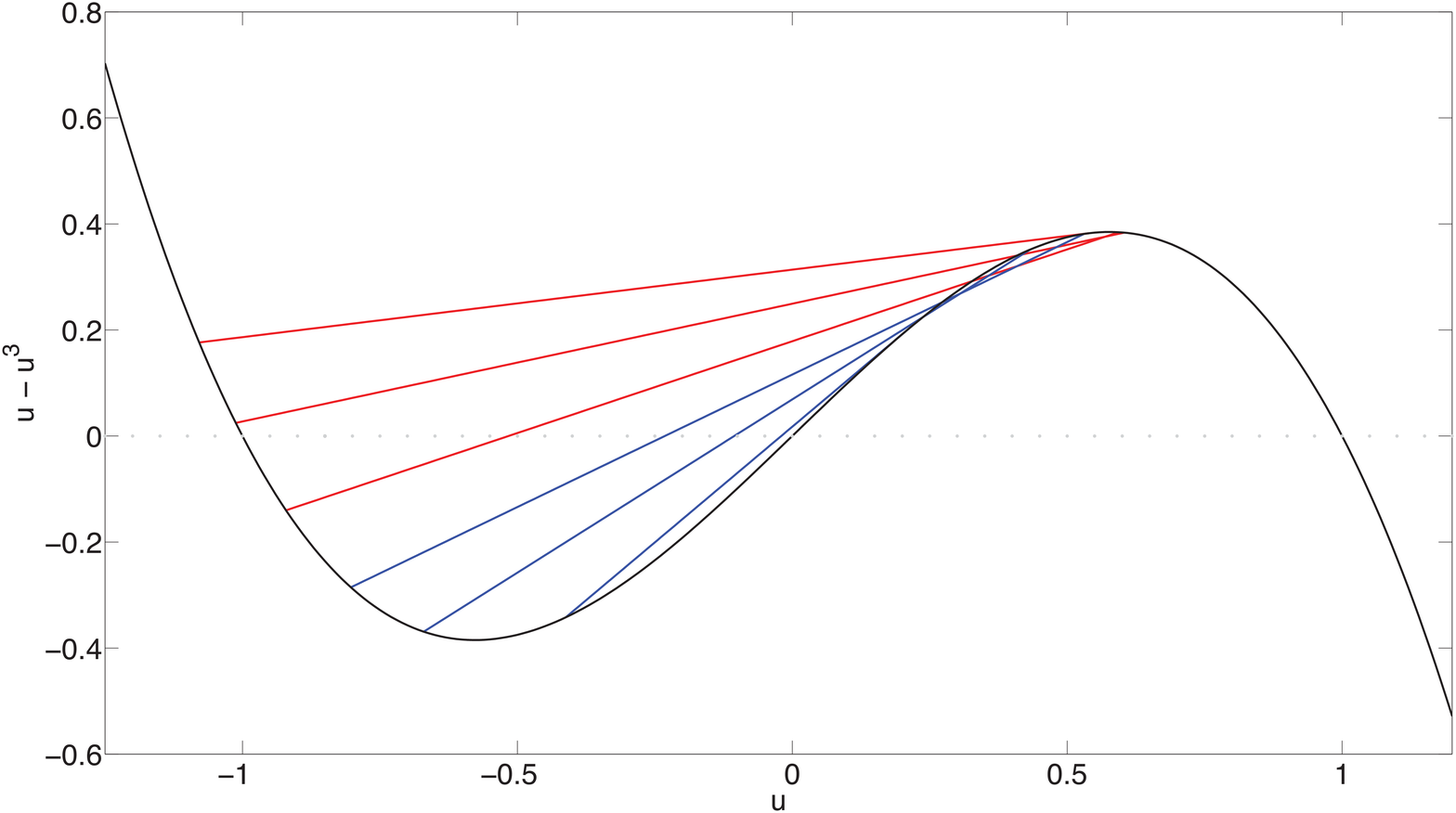}}
 \caption{ Undercompressive shocks with $\gamma=1/\sqrt{6}. $ } 
  \label{figuc1}
 \end{figure}

Let $\gamma\in (0,\sqrt{3/8})$ be fixed.  In \eqref{parametric1}, it is straightforward to check that both functions  $u_+=u_+^{(\pm)}(a)$ are monotonic in $a$ ($u_+^{(+)}$ is increasing, and $u_+^{(-)}$ is decreasing). Moreover the two functions have different ranges over the domain $\half\leq a\leq \tilde{a}: \ \ u_+^{(+)}(a)\leq u_+^{(-)}(a)$ with 
     $u_+^{(+)}(\tilde{a})=u_+^{(-)}(\tilde{a}).$ Consequently, the map $a\to u_+^{(\pm)}(a)$ is invertible, with inverse $a=\overline{a}(u_+).$ Thus, $u_-=-\overline{a}(u_+) u_+$ is uniquely defined for each $u_+$ in the union of the ranges of the pair $u_+^{(\pm)}.$ This map is not invertible for all $\gamma$ however; for some range of $\gamma,$ there are two values of $u_+$ for each $u_-$ in an interval. In Figure~\ref{figuc1}, we show chords joining points $(u_\pm, f(u_\pm))$ in the graph of  $f(u)=u-u^3,$ for a specific value of $\gamma.$ Note that each $u_+<0$ is joined to a unique $u_->0,$ but that some values of $u_-$ have two values of $u_+.$ The slope of each chord is the speed of the corresponding undercompressive wave. The slowest and fastest are tangent to the graph at $u_-,$ and indicate the range of $u_+$ and $u_-$ in the construction.

  As $a\to \half,$ the curve approaches the line $u_-=-\half u_+$ and stops. This gives upper and lower bounds on $u_+<0$   that we obtain by 
 setting $u_-=\displaystyle -\frac{u_+}{2}$ in \eqref{eq6}:
\beq
\Label{eq9}
u_+=-\sqrt{2/3\left( 1\pm\sqrt{1-8\gamma^2/3}\right)}.
\eeq

 The following theorem states that all saddle-to-saddle heteroclinic orbits are given by the above construction.
\begin{thm}\label{thm1}
 Let  $0<\gamma<\sqrt{3/8},$ and  suppose there is a saddle-saddle connection from  $u_->0$ to $u_+.$ Then   $-\sqrt{2/3\left( 1-\sqrt{1-8\gamma^2/3}\right)}>u_+>-\sqrt{2/3\left( 1+\sqrt{1-8\gamma^2/3}\right)},$ and $u_-, u_+$ are given implicitly by  equation   \eqref{parametric1}, with the value of $a$ determined uniquely from $u_+$ and $\gamma.$   
  \end{thm}
 
  {\bf Proof:} \ 
The proof relies on the following lemma, proved in \cite{shearer1995,yang}, concerning cubic vector fields of the form 
\beq \label{cubicvf}
\barr{rcl}
u^\prime&=&v\\
 v^\prime&=& b v+c(u).
\earr
\eeq
\begin{Lemma}\label{lemma1}\cite{shearer1995,yang}
Let $b\in\R,$ and suppose $c(u)$ is a cubic polynomial with three distinct zeroes $u_-, u_0, u_+,$ and such that $c'(u_\pm)>0.$
If $(u(t),v(t))$ is a solution of system \eqref{cubicvf} satisfying $ \lim_{t\to\pm\infty}(u,v)(t)=(u_\pm,0),$ respectively, then the trajectory lies on an invariant parabola for  \eqref{cubicvf}.
\end{Lemma}

  The lemma applies to system \eqref{ode1} with $b=\frac{\gamma}{\sqrt{s}}$ and $c(u)=u^3-u-(u_-^3-u_-)+s(u-u_-).$ Since the formulas \eqref{parametric1}   establish  all invariant parabolas for system \eqref{ode1}, we only need the uniqueness of $u_-$ for each $u_+,$ which we established earlier.  
 This completes the proof of the theorem. \proofend
   
\begin{figure}
   \centerline{\includegraphics[height=3.5in]{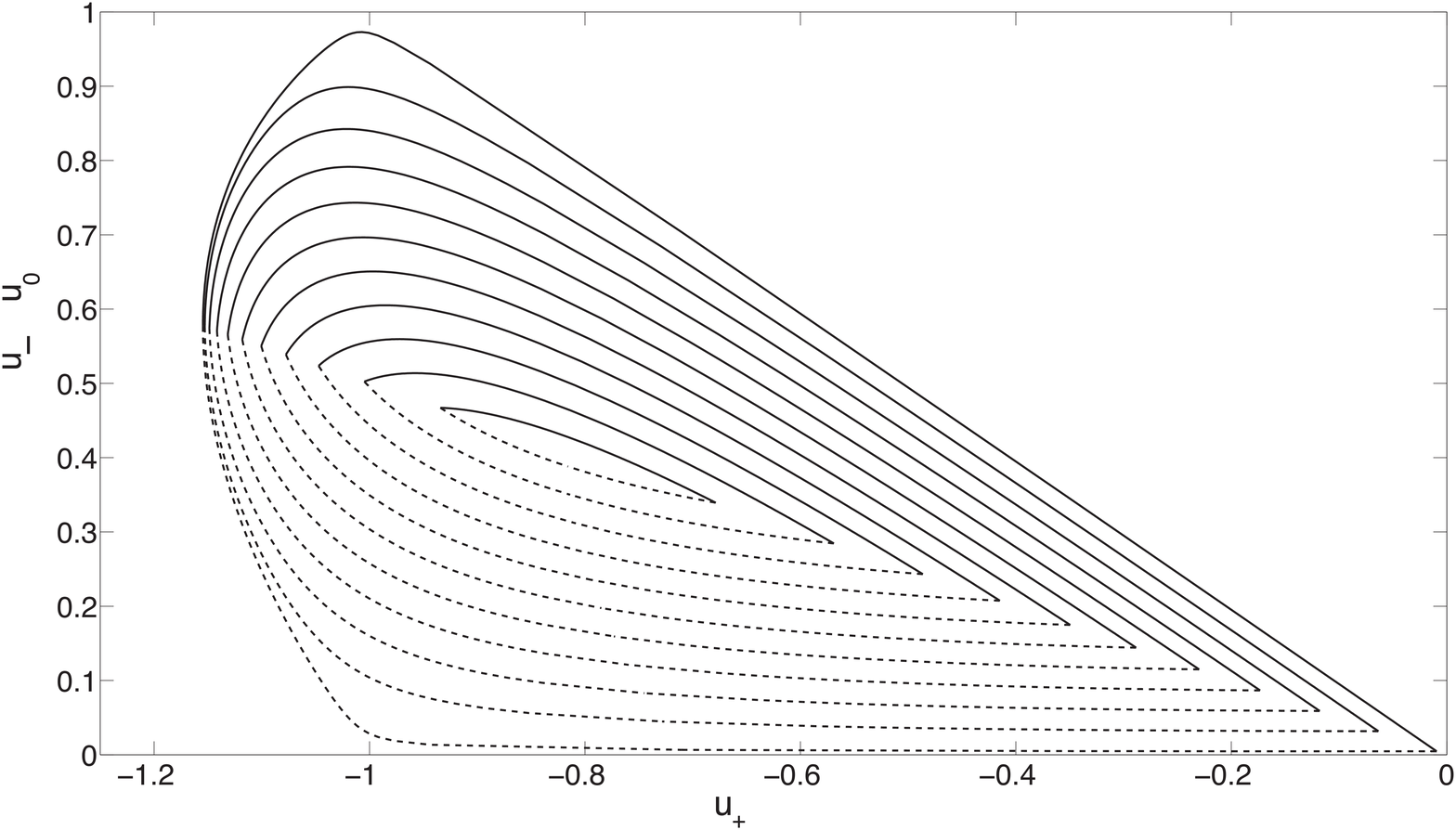}}
 \caption{$u_-$ (solid), $u_0$ (dashed) vs. $u_+$ for $\gamma=\frac{n}{10} \sqrt{\frac 38},\ n=1,...,10.$ Both signs in formula \eqref{parametric1} are used in plotting these curves.} 
  \label{figuc2}
 \end{figure}

  In Fig.~\ref{figuc2} we plot the formulas  \eqref{parametric1} together with the corresponding middle equilibria $u_0,$ given by \eqref{uzero} (with the minus sign),  for  a range of $\gamma\leq \sqrt{3/8}.$
   
\section{The Riemann Problem}\label{sec_rp}
The Riemann initial value problem for equation \eqref{cl} involves initial data with two constants $u_L, u_R:$
\beq\label{rp1}
u_t+(u-u^3)_x=0, \qquad 
u(x,0)=\begin{cases}
u_L&\text{if \ $x<0$}\\[8pt]
u_R&\text{if \ $x>0$}.
\end{cases}
\eeq
 A solution resolves the initial jump discontinuity into a combination of shocks, rarefaction waves and constants. The Riemann problem  is scale invariant, so that the solution is necessarily a function of the similarity variable $x/t$ if it is to be unique. 

\begin{figure}[ht]
\centerline{\includegraphics[height=4in]{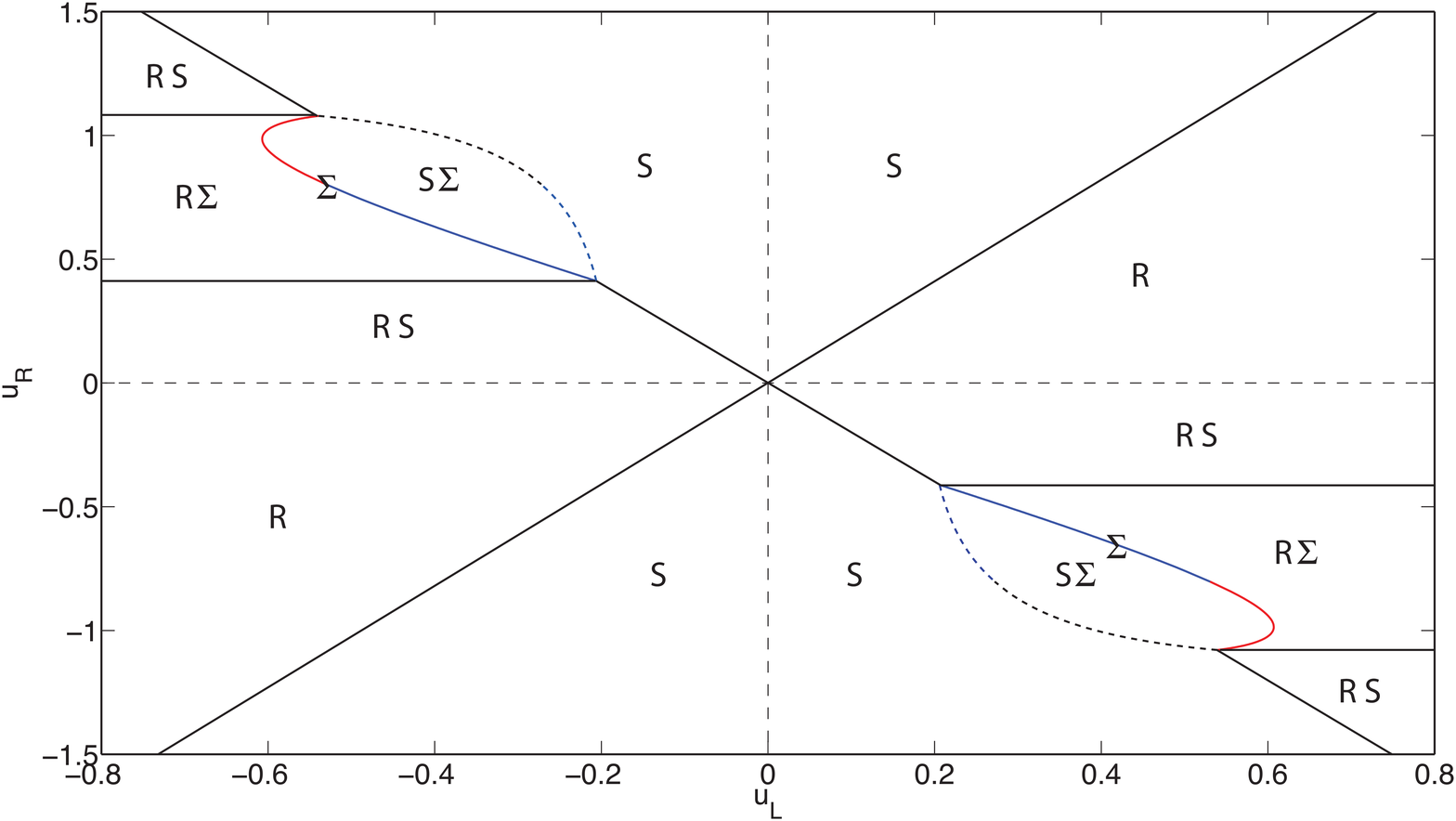}}
 \caption{   Solution of the Riemann problem, $\gamma=1/\sqrt{6}.$ \ \ R: rarefaction;  S: Lax shock;  \qquad
 $\Sigma:$ undercompressive shock.}
 \label{fig3}
 \end{figure}

Uniqueness of the solution for all initial data $u_L,u_R$ depends upon identifying a suitable condition on shock waves. The Lax entropy condition \cite{lax} requires characteristics to enter the shock from both sides, and leads to a unique solution of the Riemann problem for all initial data. However, not all the shocks satisfying the Lax condition are TW-admissible in the sense of \S\ref{sec_tw}.
 In Fig.~\ref{fig3} we show the solution of the Riemann problem for all initial conditions. The letters $R, S, \Sigma$ represent a rarefaction, Lax shock and undercompressive shock, respectively.  The figure is calculated with $\gamma=1/\sqrt{6},$ as in Fig.~\ref{figuc1}.

\begin{figure}[ht]
 \centerline{\includegraphics[height=4in]{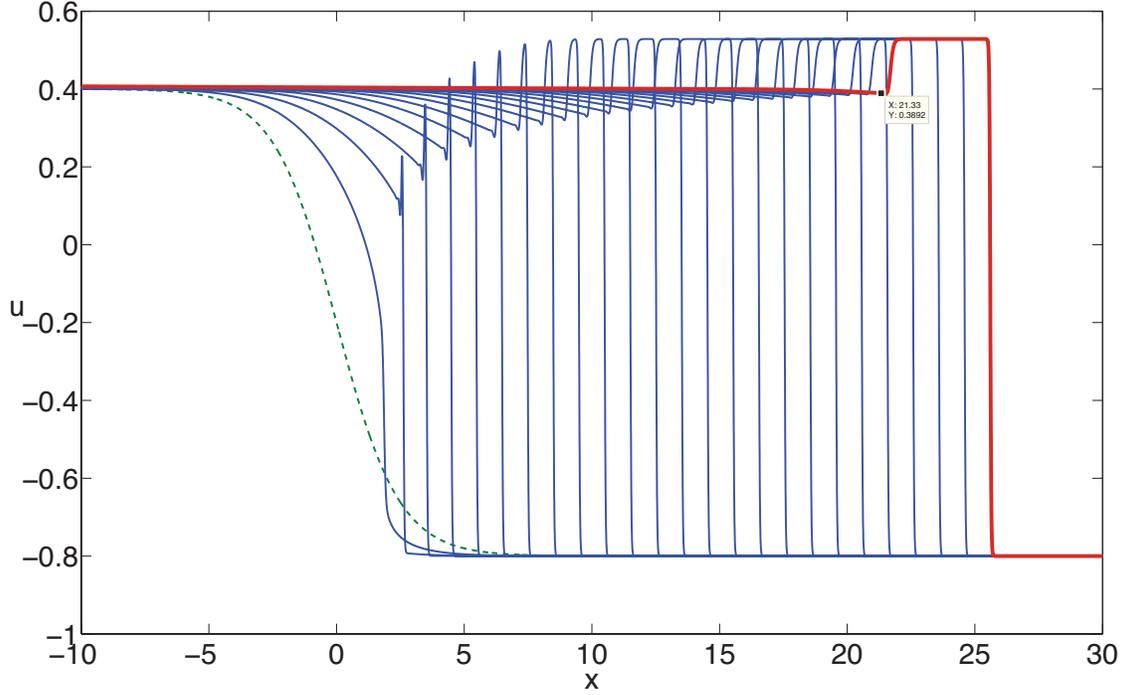}}
 \caption{   Numerical solution of the full equation \eqref{pde2} with smoothed Riemann initial data (dashed line) $u_L=0.4, u_R=-0.8;$ $\gamma=\beta/\sqrt{\mu}=1/\sqrt{6}.$ The plots are at times, $t=0,2,4,...,50.$}
 \label{fig4}
 \end{figure}

 We test the Riemann problem solution with numerical simulations in a case for which the predicted solution is an admissible Lax shock and an undercompressive shock. The initial condition  is $u(x,0)= \half\{(u_R-u_L)\tanh(\gamma x)+(u_R+u_L)\}.$
The results are shown in Fig.~\ref{fig4}. 

    \section{The p-system with Dissipation and BBM-Type Dispersion}\label{sec_system}
    Consider the system
    \begin{subequations}\Label{sysA}
    \begin{eqnarray}
    u_t-v_x&=&0\label{sysA1}\\
    v_t-\left(u^3\right)_x&=&\varepsilon v_{xx}-\varepsilon^2A u_{xxt}.\label{sysA2}
    \end{eqnarray}
    \end{subequations}

    For $\epsilon=0, $ this is the p-system, with $p(u)=-u^3.$ Characteristic speeds are $\lambda_\pm(u)=\pm\sqrt{3}u,$ so that the system loses strict hyperbolicity and genuine nonlinearity on the line $u=0, v\in \R.$ 
    The regularization on the right hand side is of the viscosity-capillarity type \cite{slemrod}, except that the capillarity term $u_{xxx}$ is replaced by the rate-dependent term $u_{xxt},$ a dispersive term analogous to the BBM dispersion in the scalar case.  
    
    There are two symmetries in the system that are useful. The obvious one is that the system is odd in the dependent variable pair $(u,v);$ the other symmetry is that the equations are invariant under  changes of sign $A\to-A, x\to -x, v\to -v.$ The final sign change can also be replaced by $u\to -u.$    
    
    We consider traveling waves of the form $(\bar{u}, \bar{v})(\xi)$ where $\xi=\frac{x-st}{\varepsilon}$. The system \eqref{sysA} then becomes (dropping the bars):
    \begin{subequations}\Label{sysB}
    \begin{eqnarray}
     -su'&=&v' \label{sysB1}\\
    -sv'-\left(u^3\right)'&=&-v''+sAu'''\label{sysB2}
    \end{eqnarray}
        \end{subequations}

    Rewriting \eqref{sysB}  as a single ODE for the single variable $u$, integrating with respect to $\xi$ and applying the boundary condition $u\to u_-$ as $\xi\to-\infty$ results in the equation
   \begin{eqnarray}
   s^2(u-u_-)-\left(u^3-u_-^3\right)=-su'+sAu''
   \end{eqnarray}
   This equation inherits the two symmetries of the PDE system \eqref{sysA}:
   \beq\label{symmetry}
   u_-\to-u_-, u\to -u, \quad \mbox{and} \quad A\to -A, \ s\to -s, \ \xi\to - \xi.
   \eeq
  The second order equation is equivalent to the system 
    \begin{subequations}\Label{sysC}
   \begin{eqnarray}
   u'&=&w\\
   sAw'&=&sw+s^2(u-u_-)-\left(u^3-u_-^3\right).
   \end{eqnarray}
       \end{subequations}
Note that the new variable $w$ is related to $v'$ in equation \eqref{sysB1}: $ w=-v'/s.$
   
  Due to the cubic nature of the system, we expect at most three equilibria, one of which is $u_-$. Equilibria have $w=0,$ and either $u=u_-,$ or 
  \begin{eqnarray}\Label{speed}
  s^2=u^2+uu_-+u_-^2.
  \end{eqnarray}
  
   We then find the other two equilibria $u_0$ and $u_+$ in terms of $u_-$ and $s$ to be
  \begin{eqnarray}\Label{quadeq}
  u_{0,+}=\frac{-u_-\pm\sqrt{u_-^2-4\left(u_-^2-s^2\right)}}{2}.
  \end{eqnarray}
  
  Notice if we consider the discriminant and substitute (\ref{speed}) for $s^2$ we find
  \begin{eqnarray}
  D=4s^2-3u_-^2=(2u+u_-)^2.
  \end{eqnarray}
 There are three real equilibria when $D>0$. The threshold $D=0$ occurs when $u_+=u_0=-\half u_-.$   When $D>0,$ the outside equilibria, which we denote $u_\pm,$ are saddle points only if $sA<0,$ as emphasized in the following Lemma:
 
 \begin{Lemma}\label{lemmasaddle}
 Suppose there are three equilibria $u_+<u_0<u_-.$ Then $u_\pm$ are saddle points if and only if $s$ and $A$ have opposite signs. 
 \end{Lemma} 
 
 {\bf Proof:} \ From \eqref{sysC} we calculate the eigenvalues at an equilibria $\tilde{u}:$
 \beq
 \lambda_\pm=\half\{A^{-1}\pm\sqrt{A^{-2}+(s^2-3\tilde{u}^2)/(sA)}\}.
 \eeq
 But for $\tilde{u}=u_\pm,$ the outside equilibria, we have $s^2<3u_\pm^2.$ Hence the result.
 \proofend
 
 A consequence of the lemma is that the only traveling waves corresponding to undercompressive shocks for system \eqref{sysA}  have speeds of only one sign, specifically opposite in sign to the sign of $A$.
 
 We are seeking a saddle-saddle connection between $u_-$ and $u_+$. 
 As in the scalar case, we seek a parabolic invariant manifold through the two equilibria: 
 \beq\label{wu1}
 w(u)=k(u-u_-)(u-u_+).
 \eeq
From the boundary conditions $u(\pm\infty)=u_\pm,$ we deduce that $w(u)=u'>0$ when $u_-<u_+,$ so that $k<0$ in that case, and $k>0$ if $u_->u_+.$ 

 Consider the equilibrium condition $s^2(u-u_-)-\left(u^3-u_-^3\right)=0$. This cubic function can be rewritten as $-(u-u_-)(u-u_+)(u-u_0)$ which is zero exactly at the equilibrium points. 
 As in \S3, we have
 \beq\label{quadeq2}
 u_0+u_++u_-=0.
 \eeq
 By the chain rule, $\frac{dw}{du}=\frac{dw}{d\xi}/\frac{du}{d\xi}=w'/u'$.  Recalling that $u'=w$,
 \begin{eqnarray}\Label{w'}
 \frac{dw}{du}=\frac{1}{sAw}\left(sw+s^2(u-u_-)-\left(u^3-u_-^3\right)\right).
 \end{eqnarray}
 Combining (\ref{w'}) with $\frac{dw}{du}=k\left(2u-(u_++u_-)\right),$ we find  
  \begin{eqnarray}
 s-\frac{1}{k}(u-u_0) =Ask(2u-(u_++u_-)).
 \end{eqnarray}
 From the constant terms,
 \begin{eqnarray}\Label{O1}
 s+\frac{1}{k}u_0=-Ask(u_++u_-)
 \end{eqnarray}
 and from the coefficient of $u,$ 
 \begin{eqnarray}\Label{Ou}
 2Ak^2=-\frac{1}{s}.
 \end{eqnarray}
 Thus, $k^2=-1/(2As),$ so that $A$ and $s$ must have the opposite signs, consistent with the assertion of Lemma~\ref{lemmasaddle}.  For definiteness, we take $A>0$ and $s<0.$ We also take the 
 positive value of $k:$ \ 
 \beq\label{k1}
 k=\frac{1}{\sqrt{-2As}}.
 \eeq 
 With these assumptions, we are seeking an invariant parabola with a trajectory from $(u_-,0)$ to $(u_+,0),$ with $u_->0$ and $w=u'<0.$ 
While the assumptions  $A>0$ and $u_->0$ may appear arbitrary, in fact the other cases with $A<0$ and/or $u_-<0$ are achieved by applying  the symmetries in \eqref{symmetry} to what follows.
 
 Multiplying (\ref{O1}) by $k$ and substituting into (\ref{Ou}) we find
 \begin{eqnarray}
 sk+u_0=\half(u_++u_-).
 \end{eqnarray}
But from (\ref{quadeq2}) we can eliminate $u_0,$ leading to
\begin{eqnarray}\label{sk2}
sk=\frac{3}{2}(u_-+u_+).
\end{eqnarray}
 Substituting for $k$ and $s$ in terms of $u_\pm,$  gives
\begin{eqnarray}
\frac{(u_+^2+u_+u_-+u_-^2)^{1/4}}{\sqrt{2A}}+\frac{3}{2}(u_++u_-)=0.
 \end{eqnarray}
  Since both terms are homogeneous in $(u_-,u_+),$ we can solve parametrically as in the scalar case. Let $u_+=bu_-.$   We then solve for $u_\pm$ as functions of $b:$
   \begin{eqnarray}\label{umsys}
  u_-=u_-(b)=\frac{2}{9(1+b)^2}\left(\frac{\sqrt{b^2+b+1}}{A}\right), \quad u_+=bu_-(b).
    \end{eqnarray}
  \begin{figure}[ht]
 \centerline{\includegraphics[height=3.5in]{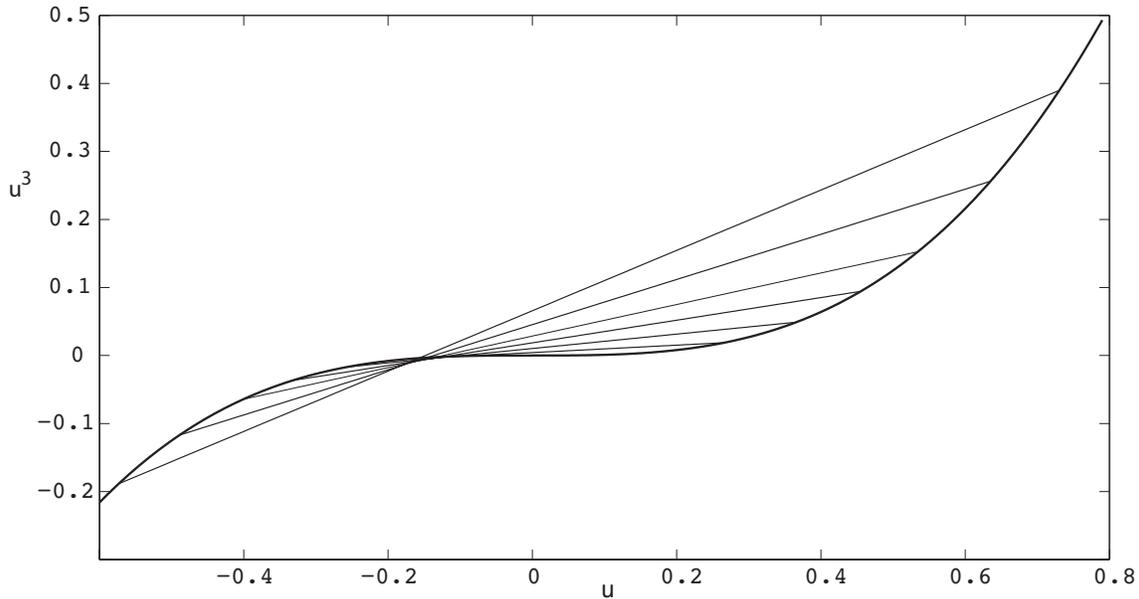}}
 \caption{ Undercompressive shocks for system \eqref{sysA} with $A=4, -0.75\leq b \leq -0.5$ in \eqref{umsys}. } 
  \label{figsysuc}
 \end{figure}

  We determine the restrictions on $b$ by examining when the number of equilibrium solutions reduces to two. We know from (\ref{quadeq}) that $u_0=-(u_++u_-)$. If $u_0=u_-$ then $u_+=-2u_-$; similarly, if $u_0=u_+$ then $u_+=\frac{-u_-}{2}$. Thus $-2u_- <u_+< \frac{-u_-}{2}$ which implies
  \begin{eqnarray}
  -2<b<-\half.
  \end{eqnarray}
  However, $u_-\to \infty$ as  $b\to -1,$ and in fact, as in the scalar case, the additional restriction $b>-1$ follows from an energy-type inequality. The only difference is that the calculation here involves $s^2$ where as in the scalar case, it is $s$ that is involved.
  Similarly, we can use Lemma~\ref{lemma1} to prove that the only saddle-saddle connections are the ones we have found in these calculations. 
  
  \begin{thm}
  Let $A>0.$ For each $u_->\frac{4\sqrt{3}}{9A},$ and each $v_-,$  there is a unique $u_+$ in the interval $-u_-<u_+<-\half u_-$ such that there is a traveling wave solution $(u,v)((x-st)/\epsilon) $ of system \eqref{sysA} satisfying    $(u,v)(\pm\infty)=(u_\pm,v_\pm),$ with speed $s=-\sqrt{u_+^2+u_+u_-+u_-^2},$ satisfying $s^2<3u_\pm^2,$ where $v_+$ is given by $v_+=v_--s(u_+-u_-).$   
     \end{thm}
     
  {\bf Proof:}  With $u_-(b)$ given by \eqref{umsys}, we have $u_-(-\half)=\frac{4\sqrt{3}}{9A}.$ Moreover,  by direct differentiation we establish easily that $u_-(b)$    is monotonically decreasing for $-1<b<-\half.$  Then Lemma~\ref{lemma1} establishes that $u_+=bu_-(b) $ is the only value of $u_+$ for which system \eqref{sysC} has an orbit from $(u,w)=(u_-,0)$  to $(u_+,0).$ The value of $v_-$ is arbitrary, since the PDE system \eqref{sysA} is invariant under translations of $v$ by a constant, and $v_+$ is given by the Rankine Hugoniot condition
  $v_+=v_--s(u_+-u_-),$ dictated by the limit of the traveling wave as $x-st\to \infty.$ This completes the proof. \proofend 
   
   \section{Discussion}
   For scalar equations and the p-system, we have introduced dispersive terms that would resemble KdV-type dispersion except that a single spatial derivative is replaced by a time derivative, in the spirit of the BBM equation \cite{benjamin1}. We find traveling wave solutions corresponding to heteroclinic orbits between saddle points. These orbits necessarily lie on invariant parabolas. Calculation of parameter values for these parabolas differs in some noticeable respects from the corresponding calculation  for the modified KdV-Burgers equation. In the scalar case, the nonlinearity is chosen carefully in order that the constant solutions are stable. This requires a balance between characteristic speeds, specifically the sign of the speeds, and the sign of the coefficient of the dispersive term. In the case considered here, with a cubic flux function, we find a bounded region of parameter values for which there can be undercompressive waves. This has implications for the Riemann problem, in that non-classical waves appear only for initial data in a restricted  region of parameter space.  Detailed properties of solutions of the Riemann problem are central to proving existence of solutions of the Cauchy problem using wave front tracking \cite{lefloch}. 
   
   In the case of systems, the situation is more complicated because genuine nonlinearity can be lost in a variety of ways. We have confined ourselves to the p-system with a homogeneous cubic function $p.$ We find that undercompressive waves can propagate only in one direction (either left or right, but not both). The direction selected depends on the sign of the dispersion coefficient. The construction used in the scalar case applies to the p-system, but the range of parameters is now unbounded. 
   
   An interesting aspect of the system case is that another natural way to incorporate a BBM-type dispersion term is to replace the term $-\varepsilon^2A u_{xxt}$ in equation \eqref{sysA} by 
   $-\varepsilon^2A v_{xxt},$ for which constant solutions are stable for $A<0,$ but are linearly unstable for high frequency perturbations if $A>0.$  However, invariant parabolas exist only for $A>0.$ Consequently, this variation in the system case requires some adjustment to the nonlinear flux function along the lines that were achieved in the scalar case by introducing $f(u)=u-u^3$ in place of the homogeneous flux $f(u)=u^3.$ 
   
 \section*{Acknowledgement}  
 Research of Michael Shearer was supported by NSF Grant DMS 0968258.

\end{document}